\documentclass[twocolumn,english]{IEEEtran}
\usepackage[T1]{fontenc}
\usepackage{babel}
\usepackage{graphicx}
\usepackage[unicode=true,
 bookmarks=true,bookmarksnumbered=true,bookmarksopen=true,bookmarksopenlevel=1,
 breaklinks=false,pdfborder={0 0 0},backref=false,colorlinks=false]
 {hyperref}
\hypersetup{pdftitle={Your Title},
 pdfauthor={Your Name},
 pdfpagelayout=OneColumn, pdfnewwindow=true, pdfstartview=XYZ, plainpages=false}

\makeatletter
\ifCLASSOPTIONcompsoc
\usepackage[caption=false,font=normalsize,labelfont=sf,textfont=sf]{subfig}
\else
\usepackage[caption=false,font=footnotesize]{subfig}
\fi

\makeatother

\begin{document}

\title{A Novel Algorithm for Linear Programming}

\author{K. Eswaran%
\thanks{Address: Dept. of Computer Science, SNIST, Jawaharlal Nehru Univ.,
Yamnampet, Ghatkesar, Hyderabad 501 301, India. Email: kumar.e @gmail.com %
}}
\maketitle
\begin{abstract}
The problem of optimizing a linear objective function, given a number
of linear constraints has been a long standing problem ever since
the times of Kantorovich, Dantzig and von Neuman. These developments
have been followed by a different approach pioneered by Khachiyan
and Karmarkar. 

In this paper we present an entirely new method for solving an old optimization
problem in a novel manner,  a technique that reduces the dimension of
the problem step by step and interestingly is recursive. A theorem
which proves the correctness of the approach is given.

The method can be extended to other types of optimization problems
in convex space, e.g. for solving a linear optimization problem subject
to nonlinear constraints in a convex region.\end{abstract}

\section{Introduction}

\IEEEPARstart{ T}{he} problem of optimizing a linear functional
subject to a set of linear constraints (the so called Linear Programming
or LP problem) has attracted many researchers; the first fundamental
contributions to the LP problem was done by Kantorovich {[}1{]} and
Dantzig {[}2{]}, who first discovered the Simplex method, which is
essentially a search method. After many years, the next break through
came through Kachiyan {[}3{]} who proved that the problem can be solved,
in theory, in polynomial time and then subsequently Karmarkar {[}4{]}
discovered a method of search involving points inside the feasible
space. Another proposed method of tackling the problem was the gravitational
method {[}5{]}. However, in spite of all these many developments {[}6,7{]}
the LP problem did not permit an easy resolution and even now a satisfactory
solution to the problem is yet to be had. The mathematician Stephen
Smale considers it as one of his 18 unsolved problems in mathematics.

In this paper\footnote{This paper was presented in: UKSIM-AMSS 14th International Conference on Modelling and Simulation,Cambridge University (Emmanuel College),U.K., March 28-30, 2012; also submitted in Arxiv.org on March 20, 2013}
we present an alternative approach to the solution of the LP problem,
the method has the advantage that an LP problem of n dimensions (excluding
{}``slack'' variables) involving an objective function of n variables,
is reduced to another LP problem in n-1 dimensions and so on. We show
and prove by a rigorous theorem, that by invoking the convexity properties
of the feasible region this stage by stage dimensional reduction of
the problem is made possible. 

The next section gives the details of the method.

\section{Description of Proposed Method}

The LP problem that will be dealt in this paper is concerned with
the task of maximizing the objective function Z defined as:

\begin{equation}
Z=d_{1}x_{1}+d_{2}x_{2}+.....d_{n}x_{n}
\end{equation}

subject to the constraints
\begin{equation}
A\, x\leq r
\end{equation}

where A is a m x n matrix and x is a n x 1 variable vector and
r is a m x 1 constant vector. Thus we are optimizing an objective
function Z, involving n variables (dimensions), given m constraint
planes.

As is well known the feasible region will be a convex polytope whose
boundaries are the constraint planes. The LP problem then consists
of searching for that point x, in the feasible region which has the
maximum value of Z. Several facts about the LP problem are well known
and bears repetition: (i) The optimum point if it exists will be a
vertex and lie on the boundary, (ii) the optimum vertex will be at
the intersection of at least n boundary planes, (iii) the feasible
region is a convex region, by which if there are any two points which
are in the feasible region then the line joining P and Q will also
lie in the feasible region and (iv) the last statement will be true
even if P and Q are on the boundary of two planes, if this happens,
then the line PQ then will either be entirely inside the feasible
region or lie on the boundary. 

In the following it is assumed that an optimum vertex to the chosen
LP problem exists. Before, we proceed further it is necessary to define
some terms: 

We will call those planes whose intersections constitute the optimum
vertex as {}``roof'' planes, from (ii) above, we can see that there
must obviously be at least n roof planes. We define $n{}_{d}$ as
the unit vector which points in the optimum direction, that is $n{}_{d}$
will have components proportional to $\left\{ d{}_{1},d_{2},d_{3},..d_{n}\right\} $.
Similarly we define $n{}_{k}$ as the unit vector which is the outward
normal of the $k^{th}$ boundary plane, in this case $n{}_{k}$ is
defined as that vector whose components are proportional to the coefficients
of the $k^{th}$row of the matrix A, hence $n{}_{k}$ is proportional
to the vector $\left\{ a{}_{k1},a_{k2},a_{k3},..a_{kn}\right\} $.
It is assumed that all the normal vectors are made to point outwards,
away from the feasible region.

Let us define the the angle $\theta_{k}$ by the dot product relationship
$ cos(\theta_{k}) = n_{k}.n_{d}$. We now define the {}``flattest
plane'' as that plane k which is such that $\theta_{k} \le \theta_{j},$ 
for $(j=1,2,...,m)$. Obviously, if k is the flattest plane then
$n_{k}.n_{d} \ge n_{j}.n_{d}$ for $(j=1,2,...,m)$

For the sake of our argument, let us assume that the optimum direction
is {}``upwards'', (there is no loss of generality in this assumption).\emph{
Now we make a crucial observation: if the problem has a single (unique)
solution, then the optimum vertex will lie on the plane which is the
{}``flattest''. }This observation follows from the fact that the
feasible region is convex, and because of this the flattest plane,\emph{
must} be one of the {}``roof'' planes, see figure. In fact, if the
flattest plane is NOT one of the roof planes, that is if a steeper
feasible plane is {}``above'' the flattest plane, then the feasible
region cannot be convex - see figure. The conclusion is true for 2-dimensions,3-dimensions
and for n-dimensions see figure. We give a formal proof in Section
V. 

From the above observation we can build an algorithm which is described:

1) We start with a properly defined LP problem such that all the constraints
are in the form given in Eq. (2), we then calculate all the outward
normals $n_{j},j=1,2,...,m $ of each of the m constraint planes. 

2) We take dot products of all the normals with the object function
direction i.e. we find all $n_{j}.n_{d},j=1,2,...,m$. By examining
each such dot product we identify the {}``flattest'' plane. If k
is the flattest plane then it will have the property $n_{k}.n_{d} \ge n_{j}.n_{d}$
for all $(j=1,2,...,m)$.

3) Now since k is the flattest plane it must contain the optimum vertex.
We then examine all the coefficients of the $k^{th}$ plane and choose
that coefficient which is the largest (say $a{}_{kn}$), we then use
the kth inequality as an equation, and get an expression for $x{}_{n}$,
in terms of $x{}_{1},$ $x{}_{2}$,....,$x{}_{n-1}$. 

4) We substitute for $x{}_{n}$using the above expression in all the
other constraint planes and delete the $k^{th}$ plane. We will have
m-1 constraint planes each of which is a function of only n-1 variables.
Next we substitute for $x{}_{n}$ in the expression for the objective
function Z, the new Z will be a linear function of only n-1 variables
$x{}_{1},$ $x{}_{2}$,....,$x{}_{n-1}$. (However, we need to retain
the equation of the deleted $k^{th}$ plane, for backsubstitution
later on). 

\emph{The idea behind elimination is simple: Since the optimum vertex
will be on the flattest plane, future searches need be conducted only
in this plane i.e. n-1 dimension space. Eliminating one variable by
using the equation of this plane ensures that a search is conducted
in this plane, but with a new objective function which does not have
this variable.}

5) The objective function and the m-1 constraints obtained in step
4 represents a LP problem of reduced dimensions. We now go back to
step 2 and reduce the problem to n-2 dimensions and so on...

6) After we have recursively reduced the problem to a single variable,
say $x{}_{1}$, we find out that value of $x{}_{1}$ which maximizes
the objective function Z which is now a function of 

this single variable and which satisfies all the constraints.

7) Having found $x{}_{1}$ the rest of the variables $x{}_{2}$,....,$x{}_{n}$
can be found by backsubstituting in the equations representing the
planes which were used for the elimination of variables, starting
from the last plane and proceeding to the first in reverse order. 

8) The value of $x{}_{1},$ $x{}_{2}$,....,$x{}_{n}$, finally obtained
represents the coordinates of the optimum vertex.

The number of steps in the reduction from an n-dimesion problem to
1 is n-1, however there is a word of caution: Every time we have finish
task 2), we must ensure that the current flattest plane is NOT a redundant
plane, by the latter we mean a plane which is entirely outside the
feasible region. In case the current flattest plane is a redundant
plane then it must be deleted and the next flattest plane should be
chosen (after testing that it is not redundant)

\section{Details of Method:}

In this section we briefly write down the various steps involved in
the method.

Denoting the objective function vector $\underline{d}$ and the vector
of the fundamental variables $\underline{x}$ as:

\begin{equation}
\underline{d}=\{d{}_{1},d_{2},d_{3},..d_{n}\},
\end{equation}

\begin{equation}
\underline{x}=\{x{}_{1},x{}_{2},x{}_{3},..x_{n}\},
\end{equation}

then the LP problem involves the task of finding out the optimum valueV
defined as the maximum value of Z, where $Z=\underline{d} . \underline{x}$
i.e

\begin{equation}
V=Max(Z)
\end{equation}

subject to the constraints :

\begin{eqnarray*}
a_{11}x_{1}+a_{12}x_{2}+a_{13}x_{3}+..+a_{1n}x_{n} & \leq & r_{1}
\end{eqnarray*}

\begin{eqnarray*}
a_{21}x_{1}+a_{22}x_{2}+a_{23}x_{3}+..+a_{2n}x_{n} & \leq & r_{2,}
\end{eqnarray*}
\begin{equation}
..................................................................................
\end{equation}

\begin{eqnarray*}
a_{m1}x_{1}+a_{m2}x_{2}+a_{m3}x_{3}+..+a_{mn}x_{n} & \leq & r_{m}
\end{eqnarray*}

Normalization: 

It will be assumed that in the above equations the objective function
vector $\underline{d}$ is a unit vector and that each row of the
constraint equations have been normalized, so that the squares of
the coefficients sum up to unity. If this has not been done (say for
the jth row) one can evaluate $N{}_{j}$ where $N_{j}= \sqrt{a_{j1}^{2}+ a_{j2}^{2}+...+a_{jn}^{2}}$
and then divide the $j^{th}$ constraint by $N{}_{j}$ and redefining
$a_{ji}/N_{j}$ by the coefficient $a_{ji}$  for each $i=1,2..,n$
and the constant $r_{j}$ by $r_{j}/N_{j}$. This procedure is adopted
for convenience because then the array $\left\{ a{}_{j1},a_{j2},a_{j3},..a_{jn}\right\} $
become the components of the normal vector $n_{j}$ to the $j^{th}$
constraint plane.

Reduction procedure:

We now demonstrate the dimension of the LP problem is reduced step
by step.

(i) Find flattest plane:

calculate $t_{r} = cos(\theta_{r})$ for all $r=(1,2,..m)$ as:

\begin{equation}
t_{r}=\sum_{i=1}^{n}a_{ri}\, d_{i}
\end{equation}

(ii) Find that plane k, such that

\begin{equation}
t_{k}\ge t_{r}\,\,\,(r=1,2,....,m)
\end{equation}

the constraint plane k then will be the flattest plane and therefore
will contain the optimum vertex. Since the optimum vertex at the intersection
of n planes, it is necessary to find the other n-1 planes chosen out
of m, whose common intersection point is the optimum vertex. Now since
we know that the optimum point lies on this plane k we enforce this
latter condition, by using the $k^{th}$constraint as an equation
and then eliminating one of the variables  $\{x{}_{1},x{}_{2},x{}_{3},..x_{n}\}$.
The variable to be eliminated will be that which has the largest coefficient
in the $k^{th}$equation - to reduce round-off errors.

At this point it is assumed, for the argument, that the kth plane
has been tested for nonredundancy, then the algorithm proceeds to
step (iii). (If the plane is redundant, see Appendix, then it must
be deleted and the next flattest but not redundant plane must be chosen).

(iii)The $k^{th}$ equation is 

\begin{equation}
a_{k1}x_{1}+a_{k2}x_{2}+a_{k3}x_{3}+..+a_{kn}x_{n}-r_{k}=0
\end{equation}

find the coeficient with the largest magnitude, say, it is the $j^{th}$
that is

\[
|a_{kj}|\geq|a_{ki}|\,\,\,(i=1,2..n)
\]

the above inequalities indicate that $a_{kj}$ is the largest coefficient
in the $k^{th}$ constraint plane , so the variable $x_{j}$ can be
eliminated from all the rest of the inequalities and objective function,
since

we can write the $k^{th}$ constraint equation as: 

\begin{equation}
x_{j}=\frac{r_{k}}{a_{kj}}-\frac{a_{k1}}{a_{kj}}\, x_{1}-\frac{a_{k2}}{a_{kj}}\, x_{2}-...\frac{a_{kn}}{a_{kj}}\, x_{n}
\end{equation}

(iv) Now we will substitute (10) for variable $x_{j}$ wherever it
occurs in all the other constraints, typically if we substitute it
in the $u^{th}$ constraint equation, the coefficients will be redefined
and it will not have the $x_{j}$. Hence, we then have (for this $u^{th}$equation):

\[
a_{u1}\leftarrow(a_{u1}-a_{uj}\,\frac{a_{k1}}{a_{kj}})
\]
\[
a_{u2}\leftarrow(a_{u2}-a_{uj}\,\frac{a_{k2}}{a_{kj}})
\]
\begin{equation}
---------\,\,\,
\end{equation}

\[
a_{uj}\leftarrow0
\]

\[
r_{u}\leftarrow(r_{u}-a_{uj}\,\frac{r_{k}}{a_{kj}})
\]

Note $a_{uj}=0$. The above replacements are done for constraints
u (u=1,2,......m), except k which is (saved in memory), but deleted
and not considered further consideration in the reduction process.

(v) The new objective function

After eliminating $x_{j}$ the $n$coefficients of the objective function
become:

\[
d_{1}\leftarrow(d_{1}-d_{j}\,\frac{a_{k1}}{a_{kj}})
\]

\[
d_{2}\leftarrow(d_{2}-d_{j}\,\frac{a_{k2}}{a_{kj}})
\]

\[
-----------
\]

\[
d_{j-1}\leftarrow(d_{j-1}-d_{j}\,\frac{a_{k(j-1)}}{a_{kj}})
\]

\[
d_{j}\leftarrow0
\]

\[
d_{j+1}\leftarrow(d_{j+1}-d_{j}\,\frac{a_{k(j+1)}}{a_{kj}})
\]

\[
-------------------
\]

\[
d_{n}\leftarrow(d_{n}-d_{j}\,\frac{a_{kn}}{a_{kj}})
\]

and we have a constant term $d_{0}$which is initially zero, becomes
nonzero
\[
d_{0}\leftarrow(d_{0} + d_{j}\,\frac{r_{k}}{a_{kj}})
\]

\begin{equation}
Z=d_{0}+d_{1}x_{1}+d_{2}x_{2}+...+d_{j-1}x_{j-1}+d_{j+1}x_{j+1}...+..d_{n}x_{n}
\end{equation}

It may be noticed that the objective function has one less varaiable:
$x_{j}$ is not present.

(vi) The equation (12) now needs to be maximized with the m-1 constraints
given by equations (11) for all planes u except k. Hence (12) and
(11) now represent an LP problem \emph{but with one less dimension
i.e. n-1.} 

Hence, we normalize the set of equations (11) and (12) as stated above.
Calculate the outward normals $n_{j}$ of the new constraint planes,
for which we need one feasible point and then we put all the constraints
in the standard form (viz.eqs (4) to (6)) and proceed to paragraph
(i) and find the new flattest plane in n-1 dimensions wrt the new
objective function.

In this manner, the problem is reduced to a single variable involving
an objective function of a single variable say $x_{1}$ which can
be maximized within the constraints.

Having, found $x_{1}$ , then the previous plane which was just eliminated,
is read from storage, this is an equation in two variables, it will
have $x_{1}$ and another variable, say $x_{2}$, using this $x_{2}$
is found. After this the plane which was eliminated second-last is
retrieved from storage, this contains an additional variable say $x_{3}$
which is then found. This back substitution process which is nothing
but a version of Gauss-Seidel elimination process, is used to find
the coordinates of the optimum vertex namely  $(x{}_{1},x{}_{2},x{}_{3},..x_{n})$.

\section{Figures}

\includegraphics[scale=0.55]{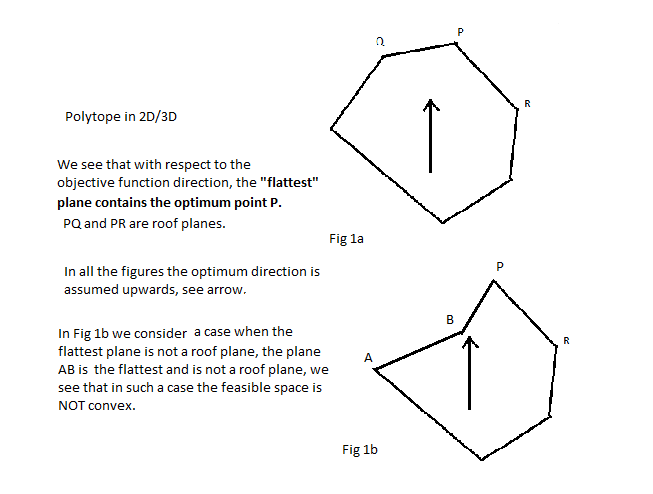}

\includegraphics[scale=0.45]{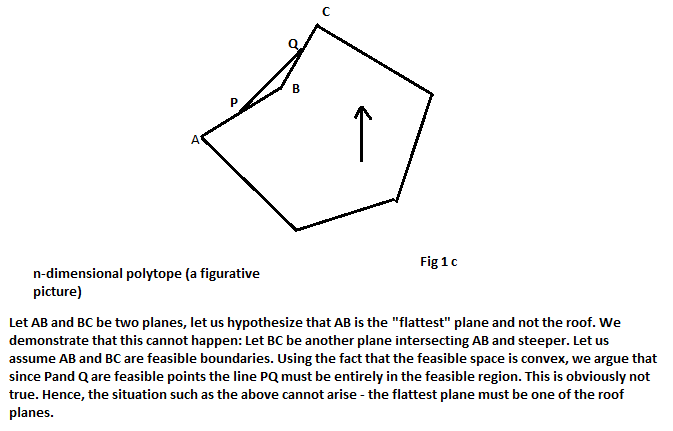}

To briefly describe the working of the algorithm, it is perhaps worthwhile
to say that as we reduce the dimensions step by step we will be choosing
the {}``flattest'' plane with respect to the present optimum direction.
For instance, suppose we are now in the r+1 th step, that is we are
now searching in n-r dimensional space, we will find that in this
space some of the planes which were having a feasible region in the
previous n-r+1 dimension space, may now become redundant, i.e. they
will lie out side the present feasible region. Of course, it is not
necessary for us to know which one unless it happens to be the {}``flattest'',
in the latter case it has to be eliminated. In other words in the
algorithm we need to test only those planes for redundancy, which
have currently qualified as {}``flattest''. Though it is acknowledged
that the test for redundancy, is a bit of a down side, the dimension
reduction along with the discovery of a roof plane at each step registers
a plus score for the algorithm. After all there are only n roof planes
to be discovered and one is being discovered in each cycle (step).

It is quite possible that the method may be extendable to other types
of optimization problems, for instance if one is to deal with the
case of optimizing a linear objective function, but under nonlinear
(but convex) constraints which still maintain a convex feasible region,
then this method is useful. Because such constraint surfaces can be
approximated by a piece-wise patchwork of planes.

\section{Theorem}

\textbf{Theorem: }The optimum vertex, if it exists lies on the flattest
plane.\medskip{}

\textbf{Proof: }

We assume that the flattest plane, in the collection of\emph{ m} planes
given in Eq. (2) is A, and given by the equation:

$ a_{k1}x_{1}+a_{k2}x_{2}+a_{k3}x_{3}+..+a_{kn}x_{n}= r_{k}$

We now show that there cannot be another plane, say B, which is {}``steeper''
than A and is a roof plane and contains a point, Q, whose optimum
value Z(Q), is higher than Z(P) of a point P, on A. Let us assume,
for the sake of argument such a B actually exists and is given by
the equation:

$ a_{j1}x_{1}+a_{j2}x_{2}+a_{j3}x_{3}+..+a_{jn}x_{n}= r_{j}$

as the proof proceeds it will be shown that such a B cannot exist.

Now since A and B are non-parallel they will {}``intersect'' in
a region R. (For the case when n=2, R is a single point, and if n=3,
R is a line and if n=4, R is a 2-d plane and so on). 

We will now perform the following coordinate transformations $ x \rightarrow x'$ :

(i) We first shift the origin to some point in the region R,

(ii) We then rotate the coordinate system such that $x'_{2}$ is along
the optimum direction, and the coordinate $x'_{1}$ is perpendicular
to it as shown, the other n-2 coordinate axis will be  $(x'_{3},x'_{4},..,x'_{n})$.
We will be able to write down the equations for A and B in terms of
a new set of coefficients such as:

Plane A:

$ sin(\alpha) x'_{1} - cos(\alpha) x'_{2}+a'_{3}x'_{3}+..+a'_{n}x'_{n}= 0$

Plane B: 

$ sin(\beta) x'_{1} - cos(\beta) x'_{2}+b'_{3}x'_{3}+..+b'_{n}x'_{n}= 0$

Note since we have chosen the new origin to be in the region R, the
constant terms in the r.h.s. of the above equations are zero; as the
origin is assumed, by our choice, to lie on both planes. 

Also since we have assumed that Plane B is steeper than plane A w.r.t.
the optimum direction $x'_{2}$, we must have:

$\beta > \alpha $

Now we show that the above inequality is untenable, with the condition
that the optimum value lies on B rather than A.

Now convexity implies that for every feasible point P which lies on
A and another feasible point Q which lies on B, the line segment PQ
must be in the feasible region. We will now show that this cannot
happen for the planes A and B chosen as above. To prove the latter
sentence all we need is to choose two feasible points P, Q lying on
A and B respectively, and show that the line PQ cannot be in the feasible
region.

We choose P as follows:

Let P be the point whose coordinate is  $(x_{P}, y_{P},0,0,...,0)$, 

that is we have chosen, $x'_{1}= x_{P},x'_{2}= y_{P},x'_{3} = 0, x'_{4} = 0,..,x'_{n} =0$.
Similarly we choose Q as the point with coordinate  $(x_{Q}, y_{Q},0,0,...,0)$,
that is we have chosen $x'_{1}= x_{Q},x'_{2}= y_{Q},x'_{3} = 0, x'_{4} = 0,..,x'_{n} =0$.
Substituting these two coordinates in the equation for their corresponding
planes we have the conditions , for P and Q to lie on A and B respectively
as:

\begin{equation}
\ensuremath{sin(\alpha)\, x_{P}-cos(\alpha)\, y_{P}=0}
\end{equation}
\begin{equation}
\ensuremath{sin(\beta)x_{Q}-cos(\beta)y_{Q}=0}
\end{equation}

In order to draw a 2-dimension figure for an n-dimensional situation,
we do as follows: consider a planar section containing the origin
O and the $ (x'_{1}, x'_{2})$ axis, in this figure, P and Q will
appear as points lying on lines OA and OB, which represent the respective
planes, A and B. Now the value of the objective function Z at points
P and Q are $Z(P) =  y_{P}$ and $Z(Q) =  y_{Q}$, since $ \beta > \alpha $,
we can see from the figure that $y_{Q} >  y_{P}$, i.e. $Z(Q) >  Z(P)$
, but the line segment PQ is outside the feasible region. Hence, plane
B, having the property as above, cannot exist, thus the theorem is
proved QED.

\includegraphics[scale=0.45]{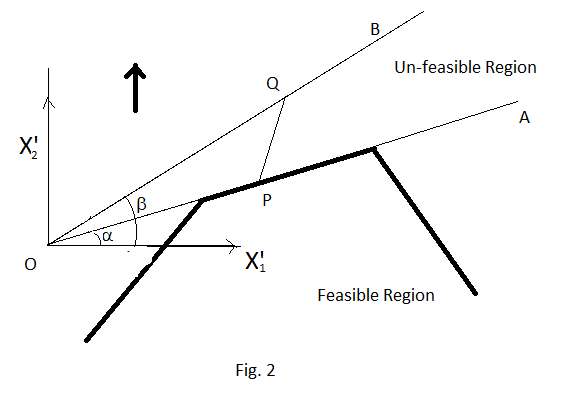}

\section{A brief on redundant planes}

In the description of the present method, we had said that the number
of steps in reducing the problem from n dimensions to 1 would be n-1,
however much depends on the actual geometry of the polytopes. It must
be ensured that the current flattest plane is not a redundant plane,
that is, it should not be a plane completely {}``above'' the feasible
region, if such a thing happens then obviously the optimum vertex
point cannot lie on it and then there is no sense using the equation
to the flattest plane for elimination of a variable/reduction - the
redundant plane must simply be deleted and the next {}``flattest''
plane must be found. The Linear Programming literature contains a
number of techniques of detecting redundant planes and these techniques
may be used.

Since we need to only test one plane at a time, the current flattest
plane to ensure that it is not redundant, perhaps, the simplest way
to begin is to use Monte Carlo techniques (e.g. see {[}8{]}). That
is randomly generate coordinates of many feasible points inside the
polytope and then from each of these points {}``draw'' straight
lines in the optimum direction, all of them will intersect some plane
or the other and out of these a few will hit the {}``flattest''
plane (see Figure 2). The points where the flattest plane is hit will
become feasible points in the n-1 dimension space upon dimension reduction
(one of the feasible points can also be used to calculate the outward
normals of the new system of planes). 

If it so happens that the flattest plane is redundant then the plane
will not be hit before another {}``lower'' plane, proving that the
plane is redundant and can be deleted. The diagram in Fig 3, shows
LM as the \textquotedblleft{}flattest\textquotedblright{} but redundant
plane, hence this will be deleted by the algorithm and the next \textquotedblleft{}flattest\textquotedblright{}
plane PQ will be retained. It is not necessary to remove all redundant
planes, but only those which happen to qualify as \textquotedblleft{}flattest\textquotedblright{}
at any stage.

This method has the advantage that it can be easily implemented by
an algorithm which is parallelizable, thus one can very efficiently
use multiple processors. 

\includegraphics[scale=0.35]{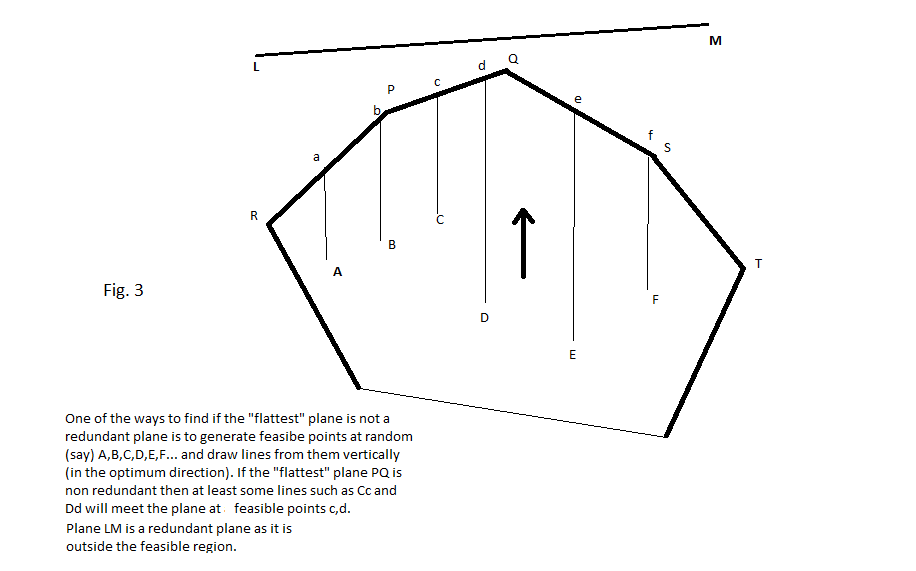}

\section{Conclusions }

A new method of solving the LP problem which is done by dimensional
reduction and is aesthetically pleasing has been found. This paper
also shows that there exists a recursive method of solving the LP
problem, a fact which was not known and therfore novel. The possible
utility of the strategy described in this paper for dealing with other
types of optimization problems involving nonlinear constraints have
also been briefly indicated.

\textbf{\Large Acknowledgments}{\Large \par}

The author thanks Professor Mukesh Eswaran for useful comments on
this paper. The help and support of the Management of SNIST in particular
Dr. K.T. Mahhe and Dr. P. Narasimha Reddy are gratefully acknowledged.

\section*{References}

(1) Kantorovich, L.V. (1939). \textquotedbl{}Mathematical Methods
of Organizing and Planning Production\textquotedbl{} Management Science,
Vol. 6, No. 4 (Jul., 1960), pp. 366\textendash{}422.

(2) G.B Dantzig, G.B.: (1947). Maximization of a linear function of
variables subject to linear inequalities. Published pp. 339\textendash{}347
in T.C. Koopmans (ed.):Activity Analysis of Production and Allocation,
New York-London 1951 (Wiley \& Chapman-Hall), also, Dantzig, George:
Linear programming and extensions. Princeton University Press and
the RAND Corporation, 1963. 

(3) Khachiyan, L.G.: (1979). \textquotedbl{}A polynomial algorithm
for linear programming\textquotedbl{} (in Russian). Doklady Akademii
Nauk SSSR 244: 1093\textendash{}1096.

(4) Karmarkar, N.: (1984). \textquotedbl{}A New Polynomial Time Algorithm
for Linear Programming\textquotedbl{}, Combinatorica, Vol 4, nr. 4,
p. 373\textendash{}395.

(5) Chang, S.Y. and Murty, K.G.: (1989) {}``The Steepest Descent
Gravitational Method for Linear Programming'', Discrete Applied Mathematics
25 211-239

(6) Murty, K.G.: (1983). Linear programming. New York: John Wiley
\& Sons Inc.. pp. xix+482. ISBN 0-471-09725-X. 

(7) Todd, M.J.: (2002). \textquotedbl{}The many facets of linear programming\textquotedbl{}.
Mathematical Programming 91 (3)

(8) Kalos, M.H. and Whitlock, Paula A.: (2008). Monte Carlo Methods.
Wiley-VCH. 
\end{document}